\documentclass{amsart}

\usepackage[T1]{fontenc}

\usepackage{graphicx}
\usepackage{url}
\usepackage{tabularx} % extra features for tabular environment
\usepackage{amsmath}  % improve math presentation
\usepackage{graphicx} % takes care of graphic including machinery
\usepackage[final]{hyperref} % adds hyper links inside the generated pdf file
\usepackage{amssymb}
\usepackage{amsmath}
\usepackage{mathtools}
\usepackage{dirtytalk}
\usepackage{stmaryrd}

\usepackage[style=numeric,sorting=none]{biblatex}
\hypersetup{
	colorlinks=true,       % false: boxed links; true: colored links
	linkcolor=blue,        % color of internal links
	citecolor=blue,        % color of links to bibliography
	filecolor=magenta,     % color of file links
	urlcolor=blue         
}
\usepackage [english]{babel}
\usepackage [autostyle, english = american]{csquotes}
\MakeOuterQuote{"}

\title{On $\mathbb{P}_\text{Berk}^1$ and Coarse-Grainings of Non-Archimedean Product Spaces}
\author{Alan Goldfarb}
\date{}

\begin{document}
\maketitle

\begin{abstract}
We determine the relationship between $\mathbb{C}_p^\times \times \mathbb{C}_p$ and a Berkovich space via an equivalence relation obtained by coarse-graining. This process also establishes a correspondence between $p$-adic fields and Bruhat-Tits trees. From this space, we give a construction of the Berkovich projective line $\mathbb{P}_\text{Berk}^1$ and outline a direction for the development of a comprehensive theory of non-Archimedean AdS/CFT and consequently, adelic AdS/CFT, enabled by this construction.
\end{abstract}

\section{Introduction}
The anti-de Sitter/Conformal Field Theory correspondence (AdS/CFT) is a conjectured relationship between the quantum theory of gravity and quantum field theory. In particular, the AdS/CFT correspondence relates a theory of gravity on the bulk of asymptotically anti-de Sitter spaces and conformal field theory on the boundary of such spaces~\cite{ramallo2013introduction}. The precise nature of this correspondence is beyond the scope of this paper, but the abstract idea is that one can translate a computation on the bulk of a surface to one on its boundary, or vice versa. One can thus make use of this correspondence to decide whether to perform a computation on the boundary of a surface or its bulk depending on which is easier. It may be the case that computation on one domain is substantially more complex than on the other, so this conjecture can simplify otherwise difficult problems. 

The AdS/CFT correspondence is classically formulated on Euclidean space, but we are often interested in discrete analogs of anti-de Sitter space. That is, those consisting of discrete geometries as opposed to continuous ones. There are a number of reasons for this, one being that integrals which are quite complicated on continuous spaces are analogous to geometric sums on discrete spaces~\cite{2017}. These sums are much easier to evaluate than their continuous counterparts, so establishing a connection between the classical AdS/CFT and one that is constructed over a discrete space may be of interest. 

Previously, the suitability of the Bruhat-Tits tree as the bulk construction in a theory of $p$-adic AdS/CFT was demonstrated in a paper entitled “$p$-adic AdS/CFT” by constructing it as a quotient of a product space of $p$-adic fields which are similar to the real fields in classical AdS/CFT corresponding to hyperbolic space~\cite{2017}. In this paper, we replicate this construction using much richer $p$-adic fields which more closely resemble the fields used in the real case. In order to construct our space of interest, we should first consider a field equipped with a non-Archimedean metric. We will begin with a fair amount of background material on $p$-adic numbers and Berkovich spaces, examine finite dimensional extensions of the $p$-adic numbers, and eventually coarse-grain $\mathbb{C}_p^\times \times \mathbb{C}_p$, where $\mathbb{C}_p$ denotes the completion of the algebraic closure of the $p$-adic numbers $\mathbb{Q}_p$, into a Berkovich space.

\section{Background}
\subsection{$p$-adic Numbers and the Bruhat-Tits Tree}
Recall that any nonzero rational number $x$ can be written uniquely in the form $x = p^m(a/b)$ where $a, b, p, m$ are integers, $a, b,$ and $p$ are relatively prime, $a > 0$, and $p$ is a prime number. Denote by $\mathbb{Q}^\times$ the set of non-zero rational numbers. For any prime number $p$, define the $p$-adic norm $|\cdot|_p$ as follows: \begin{itemize}
\item For any $x \in \mathbb{Q}^\times$, the \text{p-adic absolute value} of $x$ is 
$|x|_p = 1/p^m$ 
\item $|0|_p = 0$
\end{itemize}
For any prime number $p$, the field of $p$-adic numbers, denoted by $\mathbb{Q}_p$, is the completion of the rational numbers $\mathbb{Q}$ with respect to the $p$-adic norm. This is akin to the construction of $\mathbb{R}$ as the set of equivalence classes of Cauchy sequences of rational numbers~\cite{p-adic}. The difference is that in the $p$-adic case, we are considering sequences of rational numbers that are Cauchy with respect to the metric naturally induced by the $p$-adic norm, not the traditional absolute value. 

A non-zero $p$-adic number $x$ can be expressed as a series, called the $p$-adic expansion of $x$, as follows: 
$$x = p^{v_p(x)}\sum_{m = 0}^\infty a_m p^m$$
where each $a_m \in \mathbb{F}_p$ is a digit with the further condition that $a_0 \neq 0$. $v_p(x)$ is called the $p$-adic valuation of $x$ and is uniquely determined by the following equality: 
$$|x|_p = p^{-v_p(x)}$$
These expansions will prove to be useful as we construct discrete geometries. The set of values that $v_p$ can take on is called the valuation group of $\mathbb{Q}_p$ and is $\mathbb{Z}$~\cite{2017}.

Let us informally define the Bruhat-Tits tree for any prime number $p$ as follows. First, consider the unit sphere $\mathbb{U}_p = \{x \in \mathbb{Q}_p : |x|_p = 1 \}$. Now, consider every element $p^m\mathbb{U}_p$ for some integer $m$. The first digit is non-zero, so there are $p-1$ possibilities for it. There are then $p$ possibilities for each subsequent digit. Additionally, notice that $$\mathbb{Q}_p^\times = \bigsqcup_{m \in \mathbb{Z}} p^m\mathbb{U}_p$$
We can now regard each set $p^m\mathbb{U}_p$ as the boundary of a "bush" where each point on the bulk of a bush represents a $p$-adic number specified up to some number of digits. We may root each of these bushes in a discrete "trunk" where points on the trunk correspond to the magnitude of the elements of the bush. Therefore, there is exactly one point for each integer power of $p$ and we may require that two points on the trunk corresponding to $p^m$ and $p^n$ are next to each other if and only if $m = n+1$ or $m = n-1$. We now have a uniform tree structure with coordination number $p+1$ whose boundary is $\mathbb{Q}_p^\times$. We may adjoin $0$ and $\infty$ as the boundary points on the trunk, giving us a tree with boundary $\mathbb{P}^1(\mathbb{Q}_p) = \mathbb{Q}_p \cup \{\infty\}$. This tree structure is called the Bruhat-Tits tree and is denoted by $T_p$. Let us parameterize this tree using two coordinates $(z_0, z)$ with $z_0 = p^m$ for some integer $m$ and $z \in \mathbb{Q}_p$. $(z_0, z)$ is then the point on the branch leading up to $z$ whose digits have been specified up to the $p^{v_p(z_0)}$ place. Interestingly, $T_p$ can be realized as the quotient $PGL(2, \mathbb{Q}_p)/PGL(2, \mathbb{Z}_p)$ where $\mathbb{Z}_p = \{x \in \mathbb{Q}_p : |x|_p \leq 1 \}$ is the unit ball. In previous formulations of the $p$-adic AdS/CFT correspondence, $T_p$ has been identified with the bulk, and $\mathbb{P}^1(\mathbb{Q}_p)$ has been identified with the boundary~\cite{2017}.
\begin{figure}[htp]
    \centering
    \includegraphics[scale=1.25]{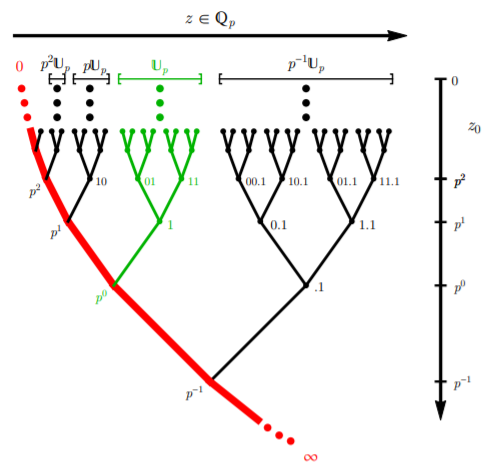}
    \caption{The Bruhat-Tits tree $T_p$ for $p=2$ with the parameterization $(z_0, z)$~\cite{2017}}
    \label{fig:1}
\end{figure}

We may construct field extensions of $\mathbb{Q}_p$ by adjoining roots of irreducible polynomials with coefficients in $\mathbb{Q}_p$. Adjoining all such roots gives the algebraic closure of $\mathbb{Q}_p$, which we denote as $\overline{\mathbb{Q}}_p$. This space is incomplete, and we denote its completion as $\mathbb{C}_p$. It is well known that for a field equipped with a particular norm, there exists a unique extension of that norm to any extension of the original field. For a field extension $K$ of $\mathbb{Q}_p$, we will use $|\cdot|_K$ to denote the extension of norm $|\cdot|_p$ to $K$~\cite{p-adic}. Note that this naturally induces a valuation function on $K$.

We will denote the degree $n$ unramified extension of $\mathbb{Q}_p$ as $\mathbb{Q}_q$ where $q = p^n$. The potential of identifying the bulk with $\mathbb{Q}_p^\times \times \mathbb{Q}_q$, a much richer space than $T_p$, has previously been studied by exploring the relationship between $\mathbb{Q}_p^\times \times \mathbb{Q}_q$ and $T_p$~\cite{2017}. This sort of identification is desirable because it may expand the class of theories one can study on the bulk by making it more closely resemble a continuum manifold while still retaining the useful tree structure. The relationship between $\mathbb{Q}_p^\times \times \mathbb{Q}_q$ and $T_p$ has been determined by performing a coarse-graining on $\mathbb{Q}_p^\times \times \mathbb{Q}_q$ to obtain an extension of the Bruhat-Tits tree: $T_q = PGL(2, \mathbb{Q}_q)/PGL(2, \mathbb{Z}_q)$~\cite{2017}. One should note that any finite extension of $\mathbb{Q}_p$ can be used to construct a corresponding extension of $T_p$.

The coarse-graining was performed using the \textit{chordal distance function} between two points $(w_0, w)$, $(z_0,z) \in \mathbb{Q}_p^\times \times \mathbb{Q}_q$ defined as follows:
$$u_p(z_0, z; w_0, w) = \frac{|(z_0-w_0, z-w)|^2_s}{|z_0w_0|_p}$$
where $|\cdot|_s$ is the supremum norm: $$|(z_0-w_0, z - w)|_s = \text{sup}\{|z_0-w_0|_p, |z-w|_q\}$$ and $|\cdot|_q$ denotes the extension of the p-adic norm to $\mathbb{Q}_q$~\cite{2017}. The main purpose of this paper is to imitate this coarse-graining to determine if a richer space has the potential to be identified with the bulk in a reformulation of the $p$-adic AdS/CFT correspondence.
\begin{figure}[htp]
    \centering
    \includegraphics[scale=1.3]{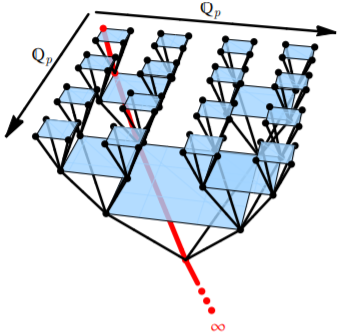}
    \caption{The Bruhat-Tits tree $T_q$ for $q=p^2$ and $p=2$~\cite{2017}.}
    \label{fig:1}
\end{figure}

\subsection{The Berkovich Projective Line}
When pursuing the objective above, we will be interested in a space which, when coarse-grained, gives a highly structured tree-like space which contains the Bruhat-Tits tree. Let us first consider a large space which contains $T_p$: the Berkovich projective line. The Berkovich projective line, denoted by $\mathbb{P}^1_{\text{Berk}}$, is an interesting space because it is analytic but also contains a tree structure. To define $\mathbb{P}^1_{\text{Berk}}$ we must first define the Berkovich affine line, $\mathbb{A}^1_\text{Berk}$, over a field $K$ which is algebraically closed and complete with respect to a non-Archimedean norm~\cite{Baker}. The closure $\overline{\mathbb{Q}}_p$ turns out not to be complete with respect to its extended norm, so we may complete $\overline{\mathbb{Q}}_p$ to form a new space, denoted by $\mathbb{C}_p$, which is both complete and algebraically closed~\cite{p-adic}. This space is constructed from $\overline{\mathbb{Q}}_p$ exactly as $\mathbb{Q}_p$ is constructed from $\mathbb{Q}$ except that the extended norm, rather than the $p$-adic norm, is used. We may then take $K = \mathbb{C}_p$ if desired. The set of points that forms $\mathbb{A}^1_\text{Berk}$ is exactly the set of all multiplicative seminorms on the polynomial ring $K[T]$. To classify these points, notice how closed disks in $K$ give rise to unique multiplicative seminorms. If $B(a, r) = \{z \in K : |z - a| \leq r\}$ is a closed disk in $K$ with $r \geq 0$, one may define a unique multiplicative seminorm $|\cdot|_{B(a, r)}$ as follows:
$$|f|_{B(a, r)} = \sup_{z \in B(a ,r)}|f(z)|$$
To obtain a compact space from this construction, it is usually necessary
to add more points. Berkovich’s classification theorem states that adding these points forms the entirety of $\mathbb{A}^1_\text{Berk}$. Denote the point in $\mathbb{A}^1_\text{Berk}$ corresponding to the seminorm $|\cdot|_x$ as $x \in \mathbb{A}^1_\text{Berk}$. The following statement of the theorem is taken directly from the referenced source: \newline

$\textbf{Berkovich’s Classification Theorem. } \textit{Every point }x \in \mathbb{A}^1_\text{Berk} \textit{corresponds to a} \newline \textit{nested sequence } B(a_1, r_1) \supseteq B(a_2, r_2) \supseteq B(a_3, r_3) \supseteq ... \textit{ of closed disks, in the sense }\newline\textit{that}$
$$|f|_x = \lim_{n\to \infty}|f|_{B(a_n,r_n)}$$
$\textit{Two such nested sequences define the same point of }\mathbb{A}^1_\text{Berk}\textit{ if and only if }$

$\textit{(a) each has a nonempty intersection, and their intersections are the same; or }$ 

$\textit{(b) both have empty intersection, and the sequences are cofinal.}$ \newline \newline
This allows us to classify points on $\mathbb{A}^1_\text{Berk}$ according to the intersection of the sequence of closed disks corresponding to points in the Berkovich affine line $B = \bigcap\limits_n B(a_n, r_n)$ as follows: \newline
\textbf{Type I:} $B$ is a point in K \newline
\textbf{Type II:} $B$ is a closed disk with radius belonging to the value group of K \newline
\textbf{Type III:} $B$ is a closed disk with radius not belonging to the value group of K \newline
\textbf{Type IV:} $B = \varnothing$ \newline \newline
Finally, we construct $\mathbb{P}^1_{\text{Berk}}$ from $\mathbb{A}^1_\text{Berk}$ by adjoining a type I point at infinity to $\mathbb{A}^1_\text{Berk}$. The point corresponding to the type II disk centered at 0 with radius 1 is called the Gauss point and is denoted by $\zeta_{Gauss}$. One should note that the sets of type I and type II points are always infinite, whereas the sets of type III and type IV points are either infinite or empty~\cite{Baker}. In our case, where the base field is $\mathbb{C}_p$, type III and IV points are present.

\begin{figure}[htp]
    \centering
    \includegraphics[scale=0.8]{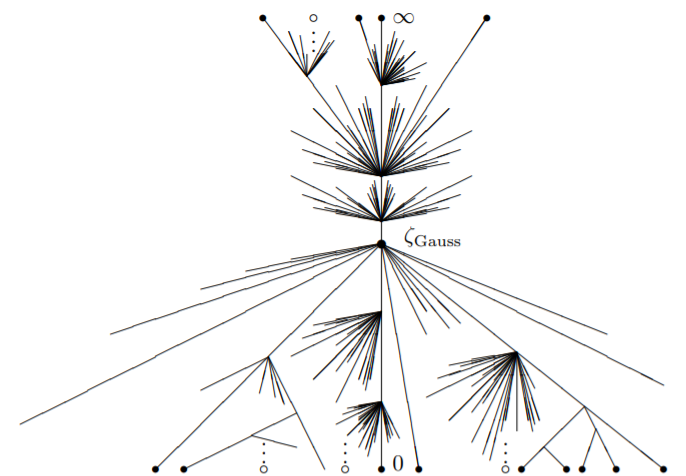}
    \caption{An artists rendition of the Berkovich Projective Line \newline $\mathbb{P}_{\text{Berk}}^1$~\cite{Baker}}
    \label{fig:1}
\end{figure}

\section{Coarse-Grainings: The Finite case}
We wish to determine the relationship between $T_p$ and product spaces of the form $K^\times \times K$ where $K$ is a field extension of $\mathbb{Q}_p$. In each case below, we extended the chordal distance function defined above on $\mathbb{Q}_p^\times \times \mathbb{Q}_q$ to the space in consideration simply by replacing the field norms used with the extended field norms associated with the corresponding extensions of $\mathbb{Q}_p$. 

To gain some intuition, we will begin with extensions of a finite degree. There are three types of finite degree extensions to consider. Each can be characterized by three quantities: the degree of the extension, its ramification index, and its residual degree. The degree of a field extension, denoted $n$, is defined as its degree as a vector space over the base field. Its ramification index, $e$, is defined as the unique natural number such that the valuation group of the field extension is $\frac{1}{e}\mathbb{Z}$. Its residual degree is simply $f = n/e$. The three types of field extensions are:
\newline \newline
i) An unramified extension: In this case, $e=1$ and $n=f$
\newline 
ii) A totally ramified extension: $e=n$ and $f=1$
\newline 
iii) A partially ramified extension: $1 < e < n$ and $n > f > 1 $~\cite{p-adic}.

\subsection{A Quadratic Extension}
We have already seen the results of coarse-graining an unramified extension in the referenced source ~\cite{2017}, so we will first consider a totally ramified extension, and work our way up to more complicated extensions using the same techniques that were used in ~\cite{2017}. Let $K = \mathbb{Q}_p(\sqrt{p})$. This is a totally ramified degree 2 extension. We may coarse-grain $K^\times \times K$ using the equivalence relation $$(z_0, z) \sim (w_0, w) \text{ iff } u_K(z_0, z; w_0, w) \leq 1$$ where $u_K$ denotes the chordal distance function on $K^\times \times K$. This function is the same as the chordal distance function on $\mathbb{Q}_p^\times \times \mathbb{Q}_q$ except that the norms $|\cdot|_p$ and $|\cdot|_q$ are replaced with the extension of $|\cdot|_p$ to $K$. We have:
$$\frac{|(z_0-w_0, z-w)|^2_s}{|z_0w_0|_K} \leq 1 \implies \left|\frac{z_0}{w_0}+\frac{w_0}{z_0}-2\right|_K \leq 1$$
It immediately follows that $\left|\frac{z_0}{w_0}\right|_K = 1$. Further, $\frac{w_0}{z_0} \in \mathbb{U}_K$ implies that $w_0 \in z_0\mathbb{U}_K$ where $\mathbb{U}_K$ is the set of units in $K$. Similarly:
$$\frac{|(z_0-w_0, z-w)|^2_s}{|z_0w_0|_K} \leq 1 \implies \left|\frac{(z-w)^2}{z_0w_0}\right|_K \leq 1 \implies \left|\frac{z-w}{z_0}\right|_K\left|\frac{z-w}{w_0}\right|_K \leq 1$$
It follows that $\left|\frac{z-w}{z_0}\right|_K \leq 1$ which implies that $w \in (z - z_0\mathbb{U}_K)$. \newline
We will denote the ball around $(z_0, z)$ containing exactly its equivalence class as $B(z_0, z)$. We have now shown that $$B(z_0, z) = (z_0\mathbb{U}_K, z - z_0\mathbb{U}_K)$$
The valuation group of $\mathbb{Q}_p(\sqrt{p})$ is $\frac{1}{2}\mathbb{Z}$, so we may choose a unique representative for the first coordinate of each equivalence class by requiring that each $z_0 = p^{\omega}$ with $\omega \in \frac{1}{2}\mathbb{Z}$. It is now evident that the first coordinate now parameterizes the trunk of the totally ramified degree two extension of the Bruhat-Tits tree. This tree is identical to the Bruhat-Tits tree $T_p$ except that the coordinate $z_0$ takes on all values of the form $p^{m/2}$ for $m \in \mathbb{Z}$ rather than just integer powers of $p$~\cite{2017}.

Let us now examine the second coordinate of each equivalence class. Consider $\pi = p^{1/2}$, a uniformizer of $K$. Recall that a uniformizer of an extended field $K$ with ramification index $e$ is an element $x \in K$ with valuation $1/e$~\cite{p-adic}. We may require that each $z$ be some element of $K$ specified up to the the $\pi^{2\omega}$ place. Formally, we require that $z \in s_{-\infty}^\omega$ where we define $s_\mu^\omega$ to be the set $\{x \in K : x = \sum_{m=\mu}^{\omega-1} a_m\pi^m, a_m \in \mathbb{F}_p\}$ when $\omega > \mu$ and $\{0\}$ otherwise~\cite{2017}. This gives that the bulk of the coarse-graining of $K^\times \times K$ is the degree 2 totally ramified extension of the Bruhat-Tits tree. We can in fact regard each equivalence class $B(z_0, z)$ as a node of this Bruhat-Tits tree.

In some cases, one may want to perform a refinement on the space resulting from the coarse-graining. The motivation for this comes from quantum gravity and is beyond the scope of this paper~\cite{2017}. To perform a refinement of this tree, let us define the set of equivalence relations $$(z_0, z) \sim_m (w_0, w) \text{ iff } u_K(z_0, z; w_0, w) \leq \pi^{-2m}$$ where $m \in \mathbb{N}$. Consider the equivalence class $B(1, 0)$ and how it splits into smaller equivalence classes under $\sim_m$. For $(z_0, z), (w_0, w) \in B(1, 0)$, we have $z_0, w_0 \in \mathbb{U}_K$. Therefore:
$$u_K(z_0, z; w_0, w) = |(z_0 - w_0, z - w)|_s^2 \leq \pi^{-2m}$$
It follows that the equivalence classes inside $B(1, 0)$ are
$$B_m(z_0, z) = (z_0 - \pi^m\mathbb{Z}_K, z - \pi^m\mathbb{Z}_K)$$
We can parameterize these equivalence classes if we require that $z_0 \in s_{0}^{m}\char`\\ s_1^m$. The set subtraction is to ensure that the first coordinate is indeed a unit in $K$. We may also require that $z \in s_{0}^{m}$. Thus, there are $|s_{0}^{m}\char`\\ s_1^m\times s_{0}^{m}| = (p^m-p^{m-1})p^m = p^{2m}(1-1/p)$ equivalence classes $B_m(z_0, z)$ in $B(1, 0)$. Because the tree structure is uniform at each node in the tree except those on the boundary, this is also the number of equivalence classes $B_m(z_0, z)$ in every other equivalence class $B(z_0, z)$. 

When $m=1$, each $B(z_0, z)$ splits into $p^2-p$ equivalence classes $B_1(z_0, z)$. In the general case, each equivalence class $B_m(z_0, z)$ splits into $\frac{p^{2(m+1)}(1-1/p)}{p^{2m}(1-1/p)} = p^2$ equivalence classes $B_{m+1}(z_0, z)$, each of which is a distance $u_p = \pi^{-2m}$ from any other equivalence class $B_{m+1}(w_0, w)$ in $B_m(z_0, z)$.

Taking the limit $m \to \infty$ of this refinement gives the original product space $K^\times \times K$. We can use an enhanced tree structure as in ~\cite{2017} to describe the refinement of the Bruhat-Tits tree into $K^\times \times K$.$^2$ First, we add $p^2-p$ edges to each node $B(z_0, z)$ on the base tree (in this case, the base tree is the degree 2 totally ramified extension of the Bruhat-Tits tree). Each of these edges corresponds to an equivalence class $B_1(z_0, z)$. From there, we add $p^2$ edges to every subsequent node $B_m(z_0, z)$ to indicate its splitting into nodes $B_{m+1}(z_0, z)$. As $m$ tends to $\infty$, the boundary of this enhanced tree structure is $K^\times \times K$ and its coordination number is $p^2+1$. Consult ~\cite{2017} for details about these refinements and enhanced tree structures.

\subsection{General Finite Extension}
We will now consider the general case, using the above coarse-graining as a motivating example. The general case is strikingly similar to this example, so little additional work is needed. Let $K$ be a degree $n$ field extension of $\mathbb{Q}_p$ with ramification index $e$ and residual degree $f = n/e$. Let us coarse-grain $K^\times \times K$. Consider the equivalence relation
$$(z_0, z) \sim (w_0, w) \text{ iff } u_K(z_0, z; w_0, w) \leq 1$$
As in the case of the quadratic extension: $B(z_0, z) = (z_0\mathbb{U}_K, z - z_0\mathbb{U}_K)$.
The valuation group of $K$ is $\frac{1}{e}\mathbb{Z}$, so we may require that each $z_0 = p^\omega$ with $\omega \in \frac{1}{e}\mathbb{Z}$.

Now consider $\pi = p^{1/e}$ a uniformizer of $K$ and let $q = p^f$. We may require that $z \in S_{-\infty}^\omega$ where we define $S_\mu^\omega$ to be the set $\{x \in K : x = \sum_{m=\mu}^{\omega-1} a_m\pi^m, a_m \in \mathbb{F}_q\}$ when $\omega > \mu$ and $\{0\}$ otherwise. This set takes on a similar numeric meaning as $s_{-\infty}^\omega$, its counterpart in the case of the quadratic extension.

The space resulting from the coarse-graining is an extension of the Bruhat-Tits tree characterized as follows: The points on the trunk are of the form $p^{m/e}$ where $m \in \mathbb{Z}$. The bulk consists of all possible specifications of all the elements of $K$. The tree therefore has coordination number $p^f+1$. The boundary of the tree consists of the points of the form $(\infty, z)$, where $z \in K$ and is therefore $\mathbb{P}^1(K)$. From this generalization, it is clear that the residual degree of $K$ determines the coordination number of the Bruhat-Tits tree resulting from its coarse-graining. The ramification index of $K$ determines the values on the trunk of the tree. This gives the following result: \newline

$\textit{Every finite field extension of } \mathbb{Q}_p \textit{ corresponds to a Bruhat-Tits tree with}\newline \textit{coordination } \textit{number } p^f+1 \textit{and whose trunk contains exactly all points of the}$ \newline \textit{form (0, $p^{m/e})$} \textit{ with }$m \in \mathbb{Z}$ \newline \newline
We may refine this tree as we did with the quadratic extension:
We will again define the following set of equivalence relations:
$$(z_0, z) \sim_m (w_0, w) \text{ iff } u_K(z_0, z; w_0, w) \leq \pi^{-2m}$$
As before, we will consider how $B(1, 0)$ splits into smaller equivalence classes under these relations. The equivalence classes inside $B(1, 0)$ take the same form as they did before. That is, $B_m(z_0 - \pi^m\mathbb{Z}_K, z - \pi^m\mathbb{Z}_K)$. We can parameterize these equivalence classes if we require that $z_0 \in S_{0}^{m}\char`\\ S_1^m$. We may also require that $z \in S_{0}^{m}$. Thus, there are $|S_{0}^{m}\char`\\ S_1^m\times S_{0}^{m}| = (q^m-q^{m-1})q^m = (p^{fm}-p^{f(m-1)})p^{fm} = p^{2fm}(1-1/p^f)$ equivalence classes $B_m(z_0, z)$ in each equivalence class $B(z_0, z)$.

When $m=1$, each $B(z_0, z)$ splits into $p^{2f}-p^{f}$ equivalence classes $B_1(z_0, z)$. In the general case, each equivalence class $B_m(z_0, z)$ splits into $\frac{p^{2f(m+1)}(1-1/p^f)}{p^{2fm}(1-1/p^f)} = p^{2f}$ equivalence classes $B_{m+1}(z_0, z)$, each of which is a distance $u_p = \pi^{-2m}$ from any other equivalence class $B_{m+1}(w_0, w)$ in $B_m(z_0, z)$.

Letting $m$ tend towards $\infty$ gives the original product space $K^\times \times K$. Again, let us use an enhanced tree structure to describe the refinement of the Bruhat-Tits tree into $K^\times \times K$: First, we add $p^{2f}-p^{f}$ edges to each node $B(z_0, z)$ on the base tree. Recall that the base tree has values $p^{l}$ for each $l \in \frac{1}{e}\mathbb{Z}$ on the trunk, and each node on the bulk is an element of $\mathbb{Q}_q$ specified up to the accuracy of the point on the trunk it emanates from. From there, we add $p^{2f}$ edges to every subsequent node $B_m(z_0, z)$ to indicate its splitting into nodes $B_{m+1}(z_0, z)$. As $m$ tends to $\infty$, the boundary of this enhanced tree structure is $K^\times \times K$ and its coordination number is $p^{2f}+1$.

\section{Constructing $\mathbb{P}_\text{Berk}^1$}
\subsection{Coarse-Graining $\mathbb{C}_p^\times \times \mathbb{C}_p$}
We begin by considering the following equivalence relation on $\overline{\mathbb{Q}}_p^\times \times \overline{\mathbb{Q}}_p$ $$(z_0, z) \sim (w_0, w) \text{ iff } u_{\overline{p}}(z_0, z; w_0, w) \leq 1$$ where $u_{\overline{p}}$ denotes the extension of $u_p$ to $\overline{\mathbb{Q}}_p^\times \times \overline{\mathbb{Q}}_p$. Let $\mathbb{U}_{\overline{p}}$ denote the unit sphere in $\overline{\mathbb{Q}}_p$ and let $\mathbb{Z}_{\overline{p}}$ denote the unit ball in $\overline{\mathbb{Q}}_p$. The equivalence relation above gives
$$B(z_0, z) = (z_0\mathbb{U}_{\overline{p}}, z - z_0\mathbb{U}_{\overline{p}})$$
The valuation group of $\overline{\mathbb{Q}}_p$ is $\mathbb{Q}$ because for any $n$, there exists some totally ramified extension of $\mathbb{Q}_p$ of degree $n$~\cite{p-adic}. This means we may require that each $z_0 = p^\omega$ with $\omega \in \mathbb{Q}$. Notice that each equivalence class will contain exactly one rational power of $p$, and that each power of $p$ is contained in an equivalence class. 

Any arbitrary element of $\overline{\mathbb{Q}}_p$ can be represented as a $p$-adic expansion by considering the smallest finite extension of $\mathbb{Q}_p$ which contains it. As with totally ramified extensions, there exists a finite unramified extension of $\mathbb{Q}_p$ with degree $n$ for every integer $n$. One may form an extension of $\mathbb{Q}_p$ of degree $n=ef$ for any ramification index $e \in \mathbb{N}$ and any residual degree $f \in \mathbb{N}$ by first taking the unramified extension of degree $f$ and then taking a degree $e$ totally ramified extension of the extended field~\cite{p-adic}. It follows that $\overline{\mathbb{Q}}_p = \{\sum_{i\in\mathbb{Q}}a_ip^i : a_i \in \overline{\mathbb{F}}_p$\} (since $\overline{\mathbb{F}}_p = \bigcup_{i=1}^\infty \mathbb{F}_{p^i}$).

Using these expansions, we notice that we can specify the second coordinate of $B(z_0, z)$ as we did before. That is, we can consider it to be an element of $\overline{\mathbb{Q}}_p$ specified up to the $p^{v_p(z_0)}$ place in the second form of the $p$-adic expansion of $z$. Formally, we may require that the second coordinate be an element of the set $\overline{S}_{-\infty}^\omega = \{\sum_{m\in M_{\omega}}a_mp^m : a_m \in \overline{\mathbb{F}}_p$\} where $M_{\omega} = \{m \in \mathbb{Q} : m < \omega\}$ and $\omega = v_p(z_0)$. This implies that the boundary of the space resulting from this coarse-graining is $\overline{\mathbb{Q}}_p$. If we now replace $\overline{\mathbb{Q}}_p$ with $\mathbb{C}_p$, the bulk of the coarse-grained space does not change: The density of $\overline{\mathbb{Q}}_p$ in $\mathbb{C}_p$ implies that every bulk point of $\mathbb{C}_p$ fits into an existing equivalence class $B(z_0, z)$ with $(z_0, z) \in \overline{\mathbb{Q}}_p$. However, the boundary of this space does change when we consider $\mathbb{C}_p$ instead of $\overline{\mathbb{Q}}_p$. We can specify an element of $K$ as $(z_0, z)$, this time with $z \in \mathbb{C}_p$. We still have that $z_0 \in \mathbb{Q}$, since the valuation groups of $\overline{\mathbb{Q}}_p$ and $\mathbb{C}_p$ are the same. One should note that unlike the parameterization described above, this is not a unique parameterization in the second coordinate. If we still require that the first coordinate be of the form $z_0 = p^\omega$, two points $(z_0, z)$ and $(w_0, w)$ are equivalent if and only if $z_0 = w_0$ and $|z-w|_{\hat{p}} < z_0^{-1}$ where $|\cdot|_{\hat{p}}$ is the extension of $|\cdot|_p$ to $\mathbb{C}_p$. Therefore, the points of the form $(\infty, z)$ are uniquely parameterized by $z \in \mathbb{C}_p$, implying that the boundary is isomorphic to $\mathbb{C}_p$.

\subsection{A Homeomorphism with $\textbf{H}^\mathbb{Q}_{\text{Berk}}$}
Let $W = \mathbb{C}_p^\times \times \mathbb{C}_p/\sim$. With the coordinates of each element of $W$  specified (although not uniquely), we may define a homeomorphism $\varphi$ between the type II points of $\mathbb{P}^1_\text{Berk}$ (denoted as $\textbf{H}^\mathbb{Q}_{\text{Berk}}$) and the bulk of $W$ as follows: Recall that a type II point of $\mathbb{P}^1_{\text{Berk}}(\mathbb{C}_p)$ is a disk $(a,r)$ centered at $a \in \mathbb{C}_p$ with radius $r \in \left|\mathbb{C}_p^\times\right|$. We may thus write $r$ as $r = p^i$ for $i \in \mathbb{Q}^\times$. Then:
$$\varphi: \textbf{H}^\mathbb{Q}_{\text{Berk}}\to \text{Bulk}(W)$$
$$ (a, p^i) \mapsto (p^{-i}, a)$$
$$ (z, p^{-\omega}) \mapsfrom (p^\omega, z)$$
if we restrict ourselves to $\overline{\mathbb{Q}}_p^\times \times \overline{\mathbb{Q}}_p$ we may send $(a, p^i)$ to $(p^{-i}, a \text{ mod }p^i)$ where $a \text{ mod }p^i$ is shorthand for the $p$-adic expansion of $a$ truncated at the $p^i$ place (so no powers of $p^j$ with $j \geq i$ are included in the $p$-adic expansion) if we wish to recover the unique parameterization of $\overline{\mathbb{Q}}_p^\times \times \overline{\mathbb{Q}}_p$ described above.

We can extend this homeomorphism to be between all of $W$ and the Type I and II points of $\mathbb{P}^1_{\text{Berk}}(\mathbb{C}_p)$ if we abuse notation so that $p^{-\infty} = 0$ and $p^\infty = \infty$. All arguments below apply to type I points using this convention if we consider type I points to be type II points with radii 0. Then the type I point corresponding to $z \in \mathbb{C}_p$ can be written in the form of a type II point as $(z, 0)$.

To verify this is a homeomorphism it suffices to check that both $\varphi$ and $\varphi^{-1}$ are well-defined and continuous. Suppose $(a_1, p^{i_1})$, $(a_2, p^{i_2})$ represent the same type II point. Then $p^{i_1} = p^{i_2}$, so we will write $p^i$ to refer to both radii. Also, it must be the case that $|a_1-a_2|_p \leq p^i$ which implies $\varphi(a_1, p^i) = \varphi(a_2, p^i)$. $\varphi^{-1}$ is trivially well-defined since each pair $(p^\omega, z)$ corresponds to an element of $W$. 

Now, recall that for any two type II points $(a_1, p^{i_1})$ and $(a_2, p^{i_2})$, $(a_1, p^{i_1})$ is on the path between $(a_2, p^{i_2})$ and the boundary of $\mathbb{P}^1_{\text{Berk}}$ if and only if $(a_1, p^{i_1}) \subset (a_2, p^{i_2})$. Note this defines a partial ordering on $\mathbb{P}^1_{\text{Berk}}$. Suppose now that $(a_1, p^{i_1}) \subset (a_2, p^{i_2})$. It follows that $(a_1, p^{i_2}) = (a_2, p^{i_2})$ so we may suppose without loss of generality that $a_1 = a_2$. Then $\varphi(a_1, p^{i_1}) = (p^{-i_1}, a_1)$ and $\varphi(a_1, p^{i_2}) = (p^{-i_2}, a_1)$. Since $(a_1, p^{i_1}) \subset (a_1, p^{i_2})$, we have that $i_2 > i_1$, so $p^{-i_2} < p^{-i_1}$. Therefore, $\varphi(a_1, p^{i_1})$ is on the path between $\varphi(a_2, p^{i_2})$ and the boundary of $W$. This argument can be applied to any two adjacent points on $\mathbb{P}_{\text{Berk}}^1$, which implies $\varphi$ is continuous. The argument for the continuity of $\varphi^{-1}$ is analogous.

It is worth noting $\varphi^{-1}$ sends the trunk to the path between 0 and infinity. In particular, $(0,0)$ maps to the type I point at 0, $(p^0,0)$ maps to the Gauss point, and $(\infty, 0)$ maps to the type I point at $\infty$. From a purely mathematical perspective, the homeomorphism established in this section is likely the most interesting result of the paper. However, we would like to construct the entirety of $\mathbb{P}^1_\text{Berk}$ in order to study $p$AdS/CFT in the future.
\subsection{Type III and IV points}
Let $W$ be the coarse-graining. Let $x = (x_0,x_1); y = (y_0,y_1) \in$ Bulk$(W)$ such that $x$ and $y$ are on the same path from the trunk to $\infty$ with $x_0 \geq y_0$. Write $x_0 = p^{i_1}$ and $y_0 = p^{i_2}$. Define a metric on Bulk$(W)$ by $\rho(x,y) = \log_p\left(\frac{p^{-i_1}}{p^{-i_2}}\right) = i_2-i_1$ (we include the intermediate step to imitate conventional notation). Let $x \lor y$ denote the highest common ancestor of $x$ and $y$ (i.e. with maximal first coordinate). In general, for two points $x, y \in K$, define $\rho(x,y) = \rho(x, x \lor y)+\rho(y,x\lor y) = i_1 + i_2 - 2i_3$. Note that the distance between adjacent points is always 1.

We can place type III points on a line segment between any type II point on the trunk, and the type I boundary. For all $w \in \mathbb{R}_{>0}\setminus\mathbb{Q}_{>0}$, $q \in \mathbb{Q}$, and $z \in \mathbb{U}_{\hat{p}}$ (where $\mathbb{U}_{\hat{p}}$ is the unit sphere in $\mathbb{C}_p$), there exists a type III point $x = (c, p^qz)$ where $c$ is the solution to the equation $w = log_p(\frac{c}{p^q}) \implies c = p^{w+q}$. The intuition is that we traverse the unique path from the point on the trunk to a boundary point $p^qz$ moving an irrational distance $w$. It is immediate that any type III point can be recovered this way since the resulting set of all type III points is then $\{(p^r, z) : r \in \mathbb{R}\setminus\mathbb{Q}, z \in \mathbb{C}_p\}$. The homeomorphism $\varphi^{-1}$ extends to a continuous bijection on these points, but it remains to show that it has a continuous inverse. Having shown this, we will have shown that adding these points makes the coarse-grained space a profinite $\mathbb{R}$-tree homeomorphic to the subset of $\mathbb{P}^1_{\text{Berk}}$ containing its type I, II, and III points.

 Let $(p^\omega, z)$ denote the point in $W$ at depth $p^\omega$ whose elements are within $p^{-\omega}$ of $z$. That is $(p^\omega, z)$ corresponds to the highest common ancestor of the boundary points in the ball of radius $p^{-\omega}$ centered at the boundary point corresponding to $z$ (this is a ball on $W$, not a ball representing a point in Berkovich space). This allows us to extend $\varphi$ to all type III points very simply (although we should recall that it is still not a unique parameterization of $W$):
$$\hat{\varphi}: \textbf{H}^\mathbb{R}_{\text{Berk}}\to \text{Bulk}(W)$$
$$ (a, p^i) \mapsto (p^{-i}, a)$$
$$ (z, p^{-\omega}) \mapsfrom (p^\omega, z)$$
We will check that this map and its inverse are well-defined and continuous. The proof that $\hat{\varphi}$ and $\hat{\varphi}^{-1}$ are well-defined is identical to the one given above for the type II points, so we will not repeat it here.

The proof of continuity is also similar. Recall the partial ordering specified on any type II or III points $(a_1, p^{i_1})$ and $(a_2, p^{i_2})$ by the condition:$(a_1, p^{i_1}) \subset (a_2, p^{i_2})$ if and only if $(a_1, p^{i_1})$ is on the path between $(a_2, p^{i_2})$ and the boundary of $\mathbb{P}^1_{\text{Berk}}$. As before, suppose that $(a_1, p^{i_1}) \subset (a_2, p^{i_2})$. It follows that $(a_1, p^{i_2}) = (a_2, p^{i_2})$ so we may suppose without loss of generality that $a_1 = a_2$. Then $\hat{\varphi}(a_1, p^{i_1}) = (p^{-i_1}, a_1)$ and $\hat{\varphi}(a_1, p^{i_2}) = (p^{-i_2}, a_1)$. Since $(a_1, p^{i_1}) \subset (a_1, p^{i_2})$, we have that $i_2 > i_1$, so $p^{-i_2} < p^{-i_1}$. Therefore, $\hat{\varphi}(a_1, p^{i_1})$ is on the path between $\hat{\varphi}(a_2, p^{i_2})$ and the boundary of $W$. This argument can be applied to any arbitrarily close points on $\mathbb{P}_{\text{Berk}}^1$, which implies $\hat{\varphi}$ is continuous given that $W$ is equipped with a suitable path metric. The argument for the continuity of $\hat{\varphi}^{-1}$ is analogous.

We have now constructed a space isomorphic to the subset of $\mathbb{P}^1_\text{Berk}$ consisting of its Type I, II, and III points. To obtain $\mathbb{P}^1_\text{Berk}$, we may simply take the unique completion of the $W$. It is already known that the completion of the type I, II, and III points of $\mathbb{P}^1_\text{Berk}$ is the entirety of $\mathbb{P}^1_\text{Berk}$, so this concludes our construction. 
\section{Future Direction}
Among the most significant drawbacks of the existing theory of $p$AdS/cFT is the difficulty of defining a spin-2 gravitational field, a necessary component in the study of quantum gravity, on the Bruhat-Tits tree. The Berkovich line is a much more refined space than the Bruhat-Tits tree and has an $\mathbb{R}$-tree structure. These properties may allow for the study of gravitational fields and dynamics on this space. The Berkovich line already has a well-studied theory of dynamics which we could apply to the physical motivations of AdS/CFT with minimal modifications~\cite{Baker}. Previously, the suitability of the Bruhat-Tits tree as the bulk construction in $p$AdS/CFT was demonstrated by presenting it as a quotient of a product space of $p$-adic fields which are similar to the real fields corresponding to hyperbolic space~\cite{2017}. We have now shown that the Berkovich projective line can be constructed using a similar quotient, but one of a much richer product space that much more closely resembles the aforementioned real fields. The next step of this project is to describe a Laplacian on the bulk of the Berkovich projective line which is appropriate for computing correlation functions. Such a Laplacian would be analogous to the real case on the type III points of the line but would be discrete on the type II points. Having done this, we will be able to define a Lagrangian and therefore calculate correlation functions. This will verify that the computations in the bulk space reproduce the conclusions of a CFT on the boundary, the first test of the correspondence on the Berkovich line. Once we have verified that the correspondence holds on this construction, a natural next step could be to seek out connections with classical AdS/CFT. In classical AdS/CFT and $p$AdS/CFT over the Bruhat-Tits tree, expressions for propagators can be given in terms of local zeta functions at the infinite and finite places, respectively, and some relations between these results are already known~\cite{2017}. A natural progression may be to derive similar expressions over the Berkovich projective line in terms of zeta functions over $\mathbb{C}_p$. This may lay a foundation for an adelic AdS/CFT.

The development of a version of non-Archimedean AdS/CFT over the Berkovich projective line would be highly impactful. The discovery of relationships between such a theory and classical AdS/CFT would give way to an adelic theory that can not only simplify calculations in the Archimedean case but also provide physical insight as to why the AdS/CFT correspondence exists to begin with. Most importantly, it would provide mathematicians and physicists the tools to properly study gravitational theories on a non-Archimedean bulk and understand their CFT analogs on the boundary.
\section*{References}
\printbibliography[heading=none]

\end{document}